\documentclass[12pt,a4paper, reqno]{amsart}
\usepackage[utf8]{inputenc}
\usepackage[T1]{fontenc}
\usepackage{amsthm}
\usepackage{amssymb}
\usepackage[abbrev]{amsrefs}
\usepackage{mathtools}
\usepackage[dvipsnames]{xcolor}
\usepackage{bm}
\usepackage[shortlabels]{enumitem}
\usepackage{hyperref}
\usepackage{bbm}

\usepackage{MnSymbol}

\newtheorem{thm}{}[section]
\newtheorem{theorem}[thm]{Theorem}
\newtheorem{corollary}[thm]{Corollary}
\newtheorem{lemma}[thm]{Lemma}
\newtheorem*{lemma*}{Lemma}
\newtheorem{proposition}[thm]{Proposition}

\theoremstyle{definition}
\newtheorem{definition}[thm]{Definition}
\theoremstyle{remark}

\newtheorem{remark}[thm]{Remark}

\numberwithin{equation}{section}
\allowdisplaybreaks

\newcommand{\FF}{\ensuremath{\mathbb{F}}}
\newcommand{\RR}{\ensuremath{\mathbb{R}}}
\newcommand{\NN}{\ensuremath{\mathbb{N}}}

\newcommand{\CC}{\ensuremath{\mathbb{C}}}

\newcommand{\xx}{\ensuremath{\bm{x}}}

\newcommand{\yy}{\ensuremath{\bm{y}}}

\newcommand{\XX}{\ensuremath{\mathbb{X}}}
\newcommand{\XB}{\ensuremath{\mathcal{X}}}
\newcommand{\VV}{\ensuremath{\mathbb{V}}}
\newcommand{\ZB}{\ensuremath{\mathcal{Z}}}
\newcommand{\YY}{\ensuremath{\mathbb{Y}}}
\newcommand{\ZZ}{\ensuremath{\mathbb{Z}}}
\newcommand{\YB}{\ensuremath{\mathcal{Y}}}

\newcommand{\Ind}{\ensuremath{\mathbbm{1}}}

\newcommand{\cG}
{\ensuremath{\mathcal{G}}}

\DeclareMathOperator{\sgn}{sign}

\DeclareMathOperator{\supp}{supp}
\DeclareMathOperator{\osc}{osc}

\newcommand{\floor}[1]{\left\lfloor #1 \right\rfloor}

\title[Quasi-greedy Markushevich bases and duality]{Quasi-greedy Markushevich bases, duality and norming subspaces}

\author[Miguel Berasategui]{Miguel Berasategui}
\address{Departamento de Matemática, Facultad de Ciencias Exactas y Naturales, Universidad de Buenos Aires, Pabellón 1, Ciudad Universitaria, (1428) Buenos Aires, Argentina}

\email{mberasategui@dm.uba.ar}

\subjclass[2020]{41A65, 46B15, 46B20.}

\keywords{thresholding greedy algorithm, quasi-greedy bases, duality, norming subspaces. }

\thanks{The author was supported by  the Grants ANPCyT PICT 2018-04104 and CONICET PIP 11220200101609CO.}

\begin{document}

\begin{abstract}
We prove that if $\XB$ is a quasi-greedy Markushevich basis of a Banach space $\XX$, its dual basis $\XB^*$ spans a norming subspace of $\XX^*$. We also prove this result for weaker forms of quasi-greediness, and study the cases of other greedy-like properties from the literature. 
\end{abstract}

\maketitle

\section{Introduction and background}
In 1999 \cite{KT1999}, Konyagin and Temlyakov introduced the Tresholding Greedy Algorithm (TGA) with respect to Schauder bases in general Banach spaces. In brief, given a basis $\XB=\left(\xx_n\right)$ 
of a Banach space $\XX$ with dual basis $\XB^*=\left(\xx_n^*\right)_{n\in \NN}\subset\XX^*$ and a vector $f\in \XX$, the TGA reorders the  set
\begin{align*}
\supp\left(f\right)=\supp\left(f\right)\left[\XB,\XX\right]:=\left\lbrace n\in \NN: \xx_n^*\left(f\right)\not=0\right\rbrace 
\end{align*}
(called the \emph{support} of $f$ with respect to $\XB$) so that in the given order (called a \emph{greedy order} of $f$) the coefficients $\xx_n^*\left(f\right)$ are non-increasing in modulus. In \cite{KT1999}, two types of Schauder bases were defined in terms of the TGA: greedy bases, and quasi-greedy ones. The former are bases in which the approximations via the greedy algorithm are the best up to a multiplicative constant, whereas the latter are those bases for which the algorithm converges, that is 
\begin{align}
&\sum_{n=1}^{\infty}\xx_{\pi\left(n\right)}^*\left(f\right)\xx_{\pi\left(n\right)}=f&&\forall f\in \XX,\label{QG1}
\end{align}
where $\pi:\NN\rightarrow \NN$ is an injective function whose range contains the support of $f$ and which chooses the coefficients in accordance to the TGA. \\
Greedy bases are strictly stronger than quasi-greedy ones and - as shown in \cite{KT1999} - the former are unconditional but the latter need not be so. Another type of basis defined in terms of the TGA - called ``almost greedy''  and which lies strictly between the two types mentioned above - was introduced by Dilworth, Kalton, Kutzarova and Temlyakov in 2003 \cite{DKKT2003}. 
These three types of bases are central in greedy approximation theory - though several other concepts of interest have been considered as well -, and have been extensively studied (see for example \cite{AA2016}, \cite{AA2017},  \cite{AABW2021}, \cite{AW2006} \cite{DKK2003}, \cite{DKOSZ2014}, \cite{DST2012} \cite{GHO2013}, \cite{Wo2000}, among others).\\
While quasi-greediness was originally defined for Schauder bases, the study of Markushevich bases with this property began shortly afterwards in \cite{Wo2000}\footnote{In \cite{Wo2000}, the concept  was also extended to quasi-Banach spaces, but we will not study that extension in this note.}, when Wojtaszczyk proved that quasi-greedy Markushevich bases are those for which the partial sums in \eqref{QG1} are uniformly bounded \cite[Theorem 1]{Wo2000} (see also \cite[Theorem 4.1]{AABW2021}, where it is shown that even the Markushevich hypothesis is not necessary and follows from a pointwise boundedness hypothesis on the sums). This mirrors the classical equivalence for Schauder bases, though with the crucial difference that for quasi-greedy Markushevich bases, the order of the sums is not fixed, but depends on each $f\in \XX$ (the order may be underdetermined by the TGA because multiple coefficients may have the same modulus, but the aforementioned equivalence holds for every order that meets the criteria \cite[Theorem 1]{Wo2000}, \cite[Theorem 2.1]{Oikhberg2018}, \cite[Theorem 4.1]{AABW2021}). \\
It is an open problem whether  every quasi-greedy Markushevich basis has a Schauder reordering (see \cite[Problem 13.6]{AABW2021}, \cite[Problem 14]{AAT2025}), and much of the research on quasi-greedy bases focuses on Schauder bases, but the study of quasi-greedy bases without a Schauder hypothesis has continued since \cite{Wo2000} (see for example \cite{AA2017}, \cite{AABW2021}, \cite{DST2012},  \cite{DKO2015}, \cite{Oikhberg2018}, among others). In this context, it is a matter of interest to study which salient properties of Schauder bases are also present in quasi-greedy Markushevich bases. \\
One well-known property of Schauder bases - which can also be useful tool for their study - is that the dual basis $\XB^*$ spans a norming subspace of $\XX^*$ (see for example \cite[Proposition 3.2.3]{AK2016}). Equivalently, if $j_{\XX}:\XX\rightarrow \XX^{**}$ is the canonical inclusion of $\XX$ in its bidual space, the restriction of the sequence $\left(j_{\XX}\left(\xx_n\right)\right)_{n\in \NN}$ to the closure of the linear span of $\XB^*$ is equivalent to $\XB$. Markushevich bases, however, do not always have this property - not even if both $\XB$ and $\XB^*$ are bounded, see for example \cite[Proposition 2.11]{AABW2021} -, and it has been so far an open question whether all quasi-greedy Markushevich bases do \cite[Problem 13.7]{AABW2021}. To our knowledge, until now the only known result in this direction was obtained by Albiac, Ansorena, Berná and Wojtaszczyk in  \cite[Theorem 2.15]{AABW2021}, where it was proved that the answer is affirmative for quasi-greedy bases that are also bidemocratic, which by \cite[Theorems 3.3 and 5.4]{DKKT2003} happens when both $\XB$ and $\XB^*$ are almost greedy. Bidemocracy was introduced in \cite{DKKT2003} to study  duality of greedy and quasi-greedy bases, and it is a strong condition under which dual bases have properties that are not found in the general case. For example, in general, the dual basis of a quasi-greedy basis does not need to be quasi-greedy, but it does when the basis is also bidemocratic (see \cite[Section 5]{DKKT2003}). As a consequence, while the result of \cite[Theorem 2.15]{AABW2021} is significant for the study of the TGA, it does not give us information on the general case of the question that concerns us. \\
The primary goal of this paper is to answer the question mentioned above, proving the norming result for all quasi-greedy Markushevich bases - and to extend it to Markushevich bases that might only have some weaker form of quasi-greediness.  Additionally, we will consider and answer similar questions for other properties of interest in greedy approximation.  \\
In the next section, we give the main definitions and notation that we need for our results or as background information, and in Section~\ref{section main results} we prove our main results. Finally, in Section~\ref{section closing remarks} we study whether other greedy-like properties from the literature guarantee norming dual bases. \\
For further background, we refer the reader to \cite{Temlyakov2008} or \cite{AABW2021}. 
\section{General setting and notation}\label{section general setting}
Throughout this paper, $\XX$ denotes an infinite dimensional, separable Banach space over a field $\FF=\RR$ or $\FF=\CC$, with dual space $\XX^*$, and $j_{\XX}: \XX\rightarrow \XX^{**}$ is the natural inclusion of $\XX$ in its bidual space. If $\ZZ$ is a subspace of $\XX^*$ and $0<c\le 1$, we say that $\ZZ$ is $c$-norming (for $\XX$) if 
\begin{align*}
&c\left\Vert f\right\Vert_{\XX}\le \sup_{\substack{f^*\in \ZZ\\ \left\Vert f^*\right\Vert=1}}\left\vert f^*\left(f\right)\right\vert &&\forall f\in \XX, 
\end{align*}
and we say $\ZZ$ is norming if there is some $0<c\le 1$ such that $\ZZ$ is $c$-norming. 

\subsection{Sets and cardinality}
Given a set $A$, 
\begin{itemize}
\item If $f,g: A\rightarrow \RR$ are functions, we write $f\lesssim g$ if there is $C>0$ such that $f(a)\le C g(a)$ for every $a\in A$. If $f\lesssim g$ and $g\lesssim f$, we write $f\approx g$.
\item If $A\subset \NN_0$, for each $B\subset \NN_0$ and $n\in \NN_0$,  $A<B$ means that $a<b$ for all $a\in A, b\in B$, and $n<A$ (resp., $A<n$) means that $\left\lbrace n\right\rbrace <A$ (resp., $A<\left\lbrace n\right\rbrace$). 
\item $\left\vert A\right\vert$ denotes the  cardinality of $A$. 
\item $A^{<\infty}$ is the set finite subsets of $A$. 
\item For each $m\in \NN$, $A^{(m)}$ and $A^{\le m}$ are the sets of subsets of $A$ of cardinality $m$ and no greater than $m$ respectively. We may use similar notation for $>$, $<$ and $\ge$. 
\item If $A\subset \RR$ and $t\in \RR$, we set $A_{\ge t}:=A\cap \left[t,+\infty\right)$; we may use similar notation for $\le$, $>$ and $<$. 
\item If $A$ is a normed space, $S_A$ denotes its unit sphere. 
\item If $A$ is a subset of a linear space, $\langle A\rangle$ denotes the linear span of $A$.
\end{itemize}
\subsection{Fundamental minimal systems}
By a \emph{fundamental minimal system} (which we may call a \emph{minimal system} or a \emph{system} when it is clear) we mean a sequence $\XB=\left(\xx_n\right)_{n\in \NN}$ whose span $\langle \XB \rangle$ is dense in $\XX$ and for which there is a sequence  $\XB^*=\left(\xx_n^*\right)_{n\in \NN}$ in the dual space $\XX^*$ that is biorthogonal to $\XB$, that is $\xx_n^*\left(\xx_k\right)=\delta_{nk}$ for each $n,k\in \NN$. We will call $\XB^*$ the \emph{dual system} of $\XB$. A fundamental minimal system $\XB$ is a \emph{Markushevich basis} for $\XX$ if $\XB^*$ is total, that is if 
\begin{align*}
&\xx_n^*\left(f\right)=0 \;\forall n\in \NN\Longrightarrow f=0, 
\end{align*}
and it is a \emph{Schauder basis} with basis constant $K_b$ if the partial sum projections are uniformly bounded and $K_b$ is the infimum $C>0$ for which 
\begin{align*}
&\left\Vert \sum_{n=1}^{m}\xx_n^*\left(f\right)\xx_n\right\Vert\le C\left\Vert f\right\Vert &&\forall m\in \NN, f\in \XX. 
\end{align*}
If, additionally, 
\begin{align*}
&\left\Vert \sum_{n=1}^{m}\varepsilon_n\xx_n^*\left(f\right)\xx_n\right\Vert\le C\left\Vert f\right\Vert &&\forall m\in \NN, f\in \XX, \left(\varepsilon_n\right)_{n\in\NN}\in S_{\FF}^{\NN},
\end{align*}
then $\XB$ is $C$-\emph{unconditional}. \\
From now on, by a \emph{basis} we mean any Markushevich basis. \\
Given a system $\XB$ for $\XX$, 
\begin{itemize}
\item The \emph{bidual system} of $\XB$ is defined as $\XB^{**}:=j_{\XX}\left(\XB\right)\big|_{\overline{\langle \XB^*\rangle}}$. 
\item For each $f\in \XX$ and $B\in \NN^{<\infty}$, $P_B\left(f\right)$ 
is the projection of $f$ on the subspace $\langle \xx_n: n\in B\rangle$, that is
\begin{align*}
P_B\left(f\right)=&\sum_{n\in B}\xx_n^*\left(f\right)\xx_n. 
\end{align*}
\item For each $f\in \XX$, we set 
\begin{align*}
\left\Vert f\right\Vert_{\infty}:=&\sup_{n\in \NN}\left\vert\xx_n^*\left(f\right)\right\vert. 
\end{align*}
\item For each $A\in \NN^{<\infty}$, $A\subset B\subset \NN$ and $\varepsilon=\left(\varepsilon_n\right)_{n\in B}\in S_{\FF}^{B}$, we define the indicator functions 
\begin{align*}
&\Ind_{\varepsilon, A}=\Ind_{\varepsilon, A}\left[\XB,\XX\right]:=\sum_{n\in A}\varepsilon_n\xx_n&&\text{and}&&&\Ind_{A}=\Ind_{A}\left[\XB,\XX\right]:=\sum_{n\in A}\xx_n.
\end{align*}
\item For each $f\in \XX$, we set $\varepsilon\left(f\right):=\left(\sgn\left(\xx_n^*\left(f\right)\right)\right)_{n\in \NN}$, where for $z\in\FF$, 
\begin{align*}
\sgn\left(z\right):=&
\begin{cases}
\frac{z}{\left\vert z\right\vert}& \text{ if }z\not=0;\\
1 & \text{ if }z=0. 
\end{cases}
\end{align*}
\item For each $m\in \NN$, we define the minimal systems $\XB_{>m}:=\left(\xx_n\right)_{n>m}$ and the spaces $\XX_{>m}:=\overline{\langle \XB_{>m}\rangle}$. We may use similar notation for $\ge$, $<$ and $\le$, and for the dual systems and spaces. 
\item If $\YB$ is a minimal system for $\YY$ and there is a bounded operator $T:\XX\rightarrow\YY$ such that $T\left(\xx_n\right)=\yy_n$ for each $n\in\NN$, we say that $\XB$ \emph{dominates} $\YB$, and we say that $\XB$ and $\YB$ are \emph{equivalent} if $T$ is an isomorphism. 
\end{itemize}
For any minimal system $\XB$, we define $\alpha_1, \alpha_2\in \left(0,\infty\right]$ by 
\begin{align*}
&\alpha_1=\alpha_1\left[\XB,\XX\right]=\sup_{n\in \NN}\left\Vert \xx_n\right\Vert;\\
&\alpha_2=\alpha_2\left[\XB^*,\XX^*\right]=\sup_{n\in \NN}\left\Vert \xx_n^*\right\Vert.
\end{align*}

\subsection{Greedy-like properties}
Greedy-like properties are usually defined in terms of greedy sums or greedy sets. We will  define $t$-greedy sets and sums for $0<t\le 1$ for greater generality (see for example \cite{AABW2021}, \cite{DKO2015}, \cite{Oikhberg2018} or \cite{Temlyakov2008}): Fix a system $\XB$ for $\XX$. Given $0<t\le 1$ and $f\in \XX$,  a set $A\subset\NN$  is \emph{$t$-greedy set} of $f$ (with respect to $\XB$) if 
\begin{align*}
&\left\vert \xx_n^*\left(f\right)\right\vert \ge t \left\vert \xx_k^*\left(f\right)\right\vert &&\forall n\in A, k\not\in A,
\end{align*}
and it is a \emph{greedy set} of $f$ if it is a $1$-greedy set. From now on, for each $m\in \NN_0$, $\cG\left(f,m,t\right)=\cG\left(f,m,t\right)\left[\XB,\XX\right]$ will denote the set of $t$-greedy sets of $f$ of cardinality $m$, whereas $\cG\left(f,t\right)=\cG\left(f,t\right)\left[\XB,\XX\right]$ will denote the set of finite $t$-greedy sets of $f$. \\
For $A\in \NN^{<\infty}$ and $f\in \XX$, $P_A\left(f\right)$ is a \emph{$t$-greedy sum} or \emph{$t$-greedy projection} of $f$ if $A\in \cG\left(f,t\right)$. \\
Let $\XB$ be a minimal system for $\XX$ and $C\ge 1$. We say that $\XB$ is
\begin{itemize}
\item  $C$-quasi-greedy if 
\begin{align*}
&\left\Vert P_A\left(f\right)\right\Vert\le C\left\Vert f\right\Vert &&\forall f\in \XX, m\in \NN, A\in \cG\left(f,m,1\right).
\end{align*}
\item $C$-almost greedy if 
\begin{align*}
&\left\Vert f-P_A\left(f\right)\right\Vert\le C\left\Vert f-P_B\left(x\right)\right\Vert &&\forall f\in \XX,  m\in \NN, A\in \cG\left(f,m,1\right), B\in \NN^{\le m}. 
\end{align*}
\item $C$-greedy if 
\begin{align*}
&\left\Vert f-P_A\left(f\right)\right\Vert\le C\left\Vert f-g\right\Vert &&\forall f\in \XX, m\in \NN,  A\in \cG\left(f,m,1\right),\\
& && g\in \langle \XB\rangle: \supp\left(g\right)\in \NN^{\le m}. 
\end{align*}
\item $C$-semi-greedy if 
\begin{align*}
&\inf_{g\in \langle \xx_n: n\in A\rangle}\left\Vert f-g\right\Vert\le C \inf_{\substack{h\in \langle\XB\rangle\\\left\vert \supp\left(h\right)\right\vert \le m}}\left\Vert f-h\right\Vert&&\forall f\in \XX, m\in \NN, A\in \cG\left(f,m,1\right). 
\end{align*}
\item $C$-bidemocratic if 
\begin{align*}
&\sup_{A\in \NN^{\le m}}\sup_{B\in \NN^{\le m}}\left\Vert\Ind_A\right\Vert \left\Vert\Ind_B\right\Vert \le C m &&\forall m\in \NN. 
\end{align*}
\item $C$-truncation quasi-greedy if 
\begin{align*}
&\left\Vert \Ind_{\varepsilon\left(f\right), A}\right\Vert \le C\left\Vert f \right\Vert &&\forall f\in \XX,  A\in \cG\left(f,1\right). 
\end{align*}
\end{itemize}
As we mentioned, greedy, almost greedy and quasi-greedy bases are central in the study of the thresholding greedy algorithm, and our focus will be on the third one of them. 
\subsection{Democracy-like properties}
Democracy was introduced in the beginning of greedy approximation theory to characterize greedy bases \cite{KT1999}, and has been extensively studied and used to connect and characterize greedy-like properties (see for instance \cite{DKKT2003}, \cite{DKK2003},  \cite{AABW2021} or \cite{AABBL2023}, among others). For example, two of the most important characterizations in greedy approximation are 
\begin{align*}
\text{greedy}\Longleftrightarrow \text{unconditional} + \text{democratic};\\
\text{almost greedy}\Longleftrightarrow \text{quasi-greedy} + \text{democratic},
\end{align*}
proved in \cite{KT1999} and \cite{DKKT2003} respectively (they were extended to quasi-Banach spaces in \cite{AABW2021}).\\
In addition to democracy, we will consider superdemocracy - a stronger property that can replace democracy in the aforementioned characterizations, see for example \cite{DKK2003} or \cite{AABW2021} -, and two weaker properties, namely conservativeness and partial democracy, introduced in \cite{DKKT2003} and \cite{AABBL2023} respectively. \\
Let $\XB$ be a minimal system for $\XX$ and $C\ge 1$. We say that $\XB$ is
\begin{itemize}
\item  $C$-\emph{democratic} if 
\begin{align*}
&\left\Vert \Ind_{B}\right\Vert\le  C\left\Vert \Ind_{A}\right\Vert && \forall A,B\in \NN^{<\infty}, |B|\le |A|. 
\end{align*}
\item  $C$-\emph{superdemocratic} if 
\begin{align*}
&\left\Vert \Ind_{\varepsilon', B}\right\Vert \le C\left\Vert \Ind_{\varepsilon, A}\right\Vert && \forall A,B\in \NN^{<\infty}: |B|\le |A|, \forall \varepsilon'\in S_{\FF}^{B},  \varepsilon\in S_{\FF}^{A}.
\end{align*}
\item $C$-\emph{conservative} if 
\begin{align*}
&\left\Vert \Ind_{B}\right\Vert \le C\left\Vert \Ind_{A}\right\Vert && \forall A,B\in \NN^{<\infty}: |B|\le |A| \text{ and }B<A. 
\end{align*}
\item $C$-\emph{partially democratic} if, for every finite set $B\subset \NN$, there is a finite set $B\subset D\subset \NN$ such that 
\begin{align*}
&\left\Vert \Ind_{B}\right\Vert\le  C\left\Vert \Ind_{A}\right\Vert && \forall A\in \NN^{<\infty}: |B|\le |A| \text{ and } A\cap D=\emptyset.
\end{align*}
\end{itemize}

\subsection{Further notation}
Throughout the paper, and in addition to the notation and definitions already introduced,  we will use standard Banach space and greedy approximation terminology (see for example \cite{AABW2021}). 
\section{Main results}\label{section main results}
To prove our results, it will be convenient to characterize the norming property that we study in a way that depends only on the system $\XB$ and the space $\XX$, without involving the dual or bidual systems or spaces. The following lemma takes care of that. 
\begin{lemma}\label{lemmasimpleeqv4}Let $\XB$ be a fundamental minimal system for $\XX$, and $K\ge 1$. Then 
\begin{enumerate}[\rm (1)]
\item \label{ida} The following are equivalent: 
\begin{enumerate}[\rm (i)]
\item \label{aleja4} For each $f\in \langle\XB\rangle$ there is $m\left(f\right)\in \NN$ such that 
\begin{align}
&\left\Vert f\right\Vert\le K \left\Vert f+g\right\Vert\label{forcont0}
\end{align}
for each  $g\in \langle \XB_{>m\left(f\right)}\rangle$. 
\item \label{funcional2}For each $f\in \langle\XB\rangle$, there is $f^*\in \langle \XB^*\rangle$ such that 
\begin{align}
&f^*\left(f\right)=\left\Vert f\right\Vert&&\text{and}&&&\left\Vert f^*\right\Vert\le K. \label{funcional3}
\end{align}
\end{enumerate}
If the above conditions hold, $\XB^*$  spans a $K^{-1}$-norming subspace. 
\item The following are equivalent:  \label{vuelta4}
\begin{enumerate}[\rm (a)]
\item \label{K-1norming} $\XB^*$ spans a $K^{-1}$ norming subspace.
\item \label{aleja5} For each $\epsilon>0$ \ref{aleja4} and \ref{funcional2} hold with $K+\epsilon$ instead of $K$ and $f\in \XX$ instead of $f\in \langle \XB\rangle$, and for some $m\left(f,\epsilon\right)$ instead of $m\left(f\right)$. 
\item \label{finitedim} For each $\epsilon>0$ and each finite dimensional subspace $\VV\subset \XX$ there is $m\left(\VV,\epsilon\right)$ such that 
\begin{align}
&\left\Vert f\right\Vert\le \left( K+\epsilon\right) \left\Vert f+g\right\Vert\label{forcont10}
\end{align}
for each $f\in \VV$ and each $g\in \overline{\langle \XB_{>m\left(\VV, \epsilon\right)}\rangle}$. 
\end{enumerate}
\end{enumerate}
\end{lemma}
\begin{proof}
\ref{ida} \ref{aleja4}$\Rightarrow$ \ref{funcional2} Given $0\not=f\in \langle \XB\rangle$, define $g^*: \langle f\rangle\oplus \langle \XB_{>m\left(f\right)}\rangle\rightarrow \FF$ by 
\begin{align*}
&g^*\left(a f+g\right)=a \left\Vert f\right\Vert&&\forall g\in \langle \XB_{>m\left(f\right)}\rangle.
\end{align*}
It follows from hypothesis that $ \left\Vert g^*\right\Vert\le K$. Let $f^*$ be a norm-preserving extension of $g^*$ to $\XX$. Then 
\begin{align*}
&\left(f^*-\sum_{n=1}^{m\left(f\right)} f^*\left(\xx_n\right)\xx_n^*\right)\left(\xx_k\right)=0&&\forall k\in \NN, 
\end{align*}
so $f^*=\sum_{n=1}^{m\left(f\right)} f^*\left(\xx_n\right)\xx_n^*$. \\
\ref{ida}\ref{funcional2} $\Rightarrow$ \ref{aleja4} Pick $f^*$ as in the statement. For any $g\in \langle \XB_{>\max\left(\supp\left(f^*\right)\right)}\rangle$ we have
\begin{align*}
 \left\Vert f\right\Vert=&f^*\left(f\right)=f^*\left(f+g\right)\le K\left\Vert f+g\right\Vert, 
\end{align*}
which completes the proof of the equivalence between \ref{aleja4} and \ref{funcional2}. \\
To complete the proof of \ref{ida}, suppose that these conditions hold, and fix  $g\in \XX$. Given $\epsilon>0$, pick $0\not=f\in \langle \XB\rangle$ with $\left\Vert f-g\right\Vert\le \epsilon$ and choose $f^*\in \langle \XB^*\rangle$ so that \eqref{funcional3} holds. Then 
\begin{align*}
\frac{\left\Vert g\right\Vert}{K}\le \frac{\epsilon+\left\Vert f\right\Vert}{K}\le \epsilon +\frac{f^*}{\left\Vert f^*\right\Vert}\left(f\right)\le 2\epsilon +\left\vert \frac{f^*}{\left\Vert f^*\right\Vert}\left(g\right)\right\vert \le 2\epsilon+\sup_{g^*\in S_{\langle \XB^*\rangle}}\left\vert g^*\left(g\right)\right\vert. 
\end{align*}
\ref{vuelta4} \ref{K-1norming} $\Rightarrow$ \ref{aleja5}  Given $0\not=f\in \XX$ and $\epsilon>0$, pick $g^*\in S_{\langle \XB^*\rangle}$ so that $g^*\left(f\right)\ge \left\Vert f\right\Vert \left(K+\epsilon\right)^{-1}$. Then \eqref{funcional3} holds for $f^*:=\frac{\left\Vert f\right\Vert g^*}{g^*\left(f\right)}$ and $K+\epsilon$ instead of $K$. This proves \ref{funcional2}, and now the same argument given before proves \ref{aleja4}, with $K+\epsilon$ instead of $K$ and for some adequate $m\left(f,\epsilon\right)$. \\
\ref{vuelta4} \ref{aleja5} $\Rightarrow$\ref{finitedim} By density, it is enough to obtain the result for $g\in \langle \XB_{>m\left(\VV, \epsilon\right)}\rangle$.  Suppose that this is false, and choose sequences $\left(f_k\right)_{k\in\NN}\subset S_{\VV}$ and $\left(g_k\right)_{k\in\NN}\subset \langle \XB\rangle$ so that for each $k$, $\supp\left(g_k\right)>k$ and \eqref{forcont10} does not hold for $f=f_k$ and $g=g_k$. We may assume that $g_k\not=-f_k$ for all $k$ and that $\left(f_k\right)_{k\in\NN}$ converges to some $f_0\in S_{\VV}$.  Choose $k_0> m\left(f_0,\frac{\epsilon}{2}\right)$.  For each $k\ge k_0$, 
\begin{align*}
K+\epsilon\le& \frac{\left\Vert f_0\right\Vert}{\left\Vert f_k+g_k\right\Vert}\le \left(K+\frac{\epsilon}{2}\right)\frac{\left\Vert f_0+g_k\right\Vert}{\left\Vert f_k+g_k\right\Vert}\le  \left(K+\frac{\epsilon}{2}\right)\left(1+\frac{\left\Vert f_0-f_k\right\Vert}{\left\Vert f_k+g_k\right\Vert}\right),
\end{align*}
which is impossible for sufficiently large $k$. \\
\ref{vuelta4} \ref{finitedim} $\Rightarrow$\ref{K-1norming}. By hypothesis and \ref{ida}, $\XB^*$ spans a $\left(K+\epsilon\right)^{-1}$-norming subspace for every $\epsilon>0$. Hence, it spans a $K^{-1}$-norming subspace. 
\end{proof}

\begin{remark}\rm \label{remarkblockbasis} It is an immediate consequence of Lemma~\ref{lemmasimpleeqv4} that if $\XB^*$ spans a $K^{-1}$-norming subspace, so does the dual system of any blocking of $\XB$. 
\end{remark}
For bases with bounded dual bases, we can weaken some of the hypotheses of  Lemma~\ref{lemmasimpleeqv4}\ref{ida}\ref{aleja4} as follows. 
\begin{corollary}\label{corollarybounded1} Let $\XB$ be a system for $\XX$, and $K\ge 1$. Suppose that $\XB^*$ is bounded and $\XB$ has the following property: For each $f\in \langle \XB\rangle$ there are $n\left(f\right)\in \NN$ and $\delta\left(f\right)>0$ such that 
\begin{align*}
&\left\Vert f\right\Vert\le K \left\Vert f+g\right\Vert\nonumber 
\end{align*}
for each $g\in \langle \XB_{>n\left(f\right)}\rangle$ with $\left\Vert g\right\Vert_{\infty}< \delta\left(f\right)$. Then $\XB$ meets the conditions of Lemma~\ref{lemmasimpleeqv4}\ref{ida}\ref{aleja4}. Thus, $\XB^*$ spans a $K^{-1}$ norming subspace. 
\end{corollary}
\begin{proof}
Suppose that the statement is false, pick $f_0\in \langle \XB\rangle$ so that the conditions of Lemma~\ref{lemmasimpleeqv4}\ref{ida}\ref{aleja4} do not hold for $f=f_0$, and let $n_1:=n\left(f_0\right)+\max\left(\supp\left(f_0\right)\right)$. By assumption, there is $g_1\in \langle \XB_{n> n_1}\rangle$ such that \eqref{forcont0} fails for $g=g_1$. Let $n_2:=\max\left(\supp\left(g_1\right)\right)$. As before, there is $g_2\in \langle \XB_{n> n_2}\rangle$ such that \eqref{forcont0} fails for $g=g_2$. Inductively, we construct a sequence $\left(g_j\right)_{j\in \NN}\subset \langle \XB\rangle$ such that for each $j\in \NN$, 
\begin{align*}
\supp\left(f_0\right)\cup\left\lbrace n\left(f_0\right)\right\rbrace<\supp\left(g_j\right)<\supp\left(g_{j+1}\right),
\end{align*}
 and \eqref{forcont0} does not hold for $g=g_j$. For each $n,j\in \NN$,
\begin{align}
\left\vert \xx_n^*\left(g_j\right)\right\vert\le& \left\vert \xx_n^*\left(f_0+g_j\right)\right\vert\le \left\Vert \xx_n^*\right\Vert \left\Vert f_0+g_j\right\Vert <\frac{\alpha_2\left\Vert f_0\right\Vert}{K}.\nonumber
\end{align}
Let 
\begin{align*}
&l_0:=1+\floor{\frac{\alpha_2 \left\Vert f_0\right\Vert}{\delta\left(f_0\right) K}}&&\text{and}&&& g_0:=\sum_{j=1}^{l_0}\frac{g_j}{l_0}.
\end{align*}
As $g_0\in \langle \XB_{n>n\left(f_0\right)}\rangle$ and $\left\Vert g_0\right\Vert_{\infty}< \delta\left(f_0\right)$, we have
\begin{align*}
\left\Vert f_0\right\Vert\le& K\left\Vert f_0+g_0\right\Vert=K\left\Vert \sum_{j=1}^{l_0}\frac{f_0+g_j}{l_0}\right\Vert \le K\sum_{j=1}^{l_0}\frac{1}{l_0}\left\Vert f_0+g_j\right\Vert <\left\Vert f_0\right\Vert,
\end{align*}
and the proof is complete. 
\end{proof}
Next, we prove our first theorem. While the most important objects of study in this context are quasi-greedy Markushevich bases, we are also interested in systems with weaker properties that still involve the uniform boundedness or the convergence of some greedy or $t$-greedy projections - some of which have been studied in the literature. For example, it is known that the normalized Haar system in $L_1\left[0,1\right]$ is not quasi-greedy \cite[Remark 6.3]{DKWo2002}, but Gogyan \cite{Gogyan2009} proved that for each $f\in \XX$ and each $0<t<1$, there is a sequence of projections $\left(P_{A_j\left(t,f\right)}\left(f\right)\right)_{j\in \NN}$ with $A_j\in \cG\left(f, j, t \right)$ for each $j\in \NN$ that is uniformly bounded and converges to $f$ \cite[Proposition 3.1, Theorem 3.1]{Gogyan2009}. Later in  \cite{DKSWo2012}, Dilworth et al. generalized the idea and studied (Schauder) bases for which there is uniform boundedness or convergence of some $t$-greedy sums. In a different context,  in \cite{Oikhberg2018} Oikhberg introduced and studied Markushevich bases for which there is uniform boundedness and convergence of all of the greedy sums (or $t$-greedy sums, for some fixed $0<t<1$) with cardinality restricted to a fixed proper infinite subset of $\NN$; some of these bases fail to be quasi-greedy \cite[Proposition 3.1 and 3.2]{Oikhberg2018} (see also \cite[Propositions 6.8 and 6.9]{BB2024}). \\
Our next theorem shows that if a minimal system has any of the above mentioned uniform boundedness properties - or even less -, then its dual system spans a norming subspace, with constant that depends only on the given uniform boundedness constant. 
\begin{theorem}\label{theoremQGequivbidual}Let $\XB$ be a system for $\XX$ with dual system $\XB^*$, $\left(d_j\right)_{j\in \NN}$ an unbounded sequence of positive integers, $\left(t_j\right)_{j\in\NN}$ a scalar sequence with $0<t_j<1$ for all $j\in \NN$, and $K\ge 1$. Suppose that for every $f\in \langle \XB\rangle$ and every $j\in \NN$ there is $A\in \cG\left(f,d_j,t_j\right)$ such that $\left\Vert P_A\left(f\right)\right\Vert \le K\left\Vert f\right\Vert$. Then $\XB^*$ spans a $K^{-1}$-norming subspace.  \\
Moreover, if $d_j=j$ for all $j\in \NN$, then for each $f\in \langle \XB\rangle$ there is $f^*\in \langle \XB^*\rangle$ such that 
\begin{align*}
&f^{*}\left(f\right)=\left\Vert f\right\Vert&&\text{and}&&&\left\Vert f^*\right\Vert\le K. 
\end{align*}
In particular, this holds if $\XB$ is $K$-quasi-greedy,
\end{theorem}
\begin{proof}
First note that for each $f\in \NN$, $\cG\left(f,d_1,t_1\right)\not=\emptyset$, so $\left(\xx_{n}^*\left(f\right)\right)_{n\in\NN}$ is bounded. Thus, by the uniform boundedness principle it follows that $\XB^*$ is bounded. \\
Now, given $\epsilon>0$, pick $0\not=f\in \langle \XB\rangle$ and $j_0\in \NN$ so that $d_{j_0}>\left\vert\supp\left(f\right)\right\vert$, then choose  $A>\supp\left(f\right)$ so that $\left\vert \supp\left(f\right)\cup A\right\vert=d_{j_0}$, and then take $0< \eta<\mu<1$ so that
\begin{align*}
&\frac{K}{1-\mu}\le K+\epsilon&&\text{and}&&&  K \eta \left\Vert \Ind_A\right\Vert \le \mu\left\Vert f\right\Vert.
\end{align*}
Now pick $\varepsilon\in \left\lbrace -1,1\right\rbrace$ so that 
\begin{align*}
&\left\Vert f\right\Vert\le \left\Vert f+\varepsilon \eta \Ind_{A}\right\Vert,
\end{align*}
set $f_0:=f+\varepsilon \eta\Ind_A$, and let
\begin{align*}
&\delta\left(f_0\right):=t_{j_0}\min_{n\in \supp\left(f_0\right)}\left\vert \xx_n^*\left(f_0\right)\right\vert &&\text{and}&&&n\left(f_0\right):=\max\left(\supp\left(f_0\right)\right). 
\end{align*}
Note that if $g\in \langle \XB_{>n\left(f_0\right)}\rangle$ and $\left\Vert g\right\Vert_{\infty}< \delta\left(f_0\right)$, then \begin{align*}
\cG\left(f_0+g, d_{j_0}, t_{j_0}\right)=\left\lbrace \supp\left(f_0\right)\right\rbrace,
\end{align*} 
because $\cG\left(f_0+g, d_{j_0}, t_{j_0}\right)\not=\emptyset$ and, if $k\in \supp\left(f_0+g\right)\setminus\supp\left(f_0\right)$, then $k>n\left(f_0\right)$ so 
\begin{align*}
\left\vert \xx_k^*\left(f_0+g\right)\right\vert\le& \left\Vert g\right\Vert_{\infty}< \delta\left(f_0\right)=t_{j_0}\min_{n\in \supp\left(f_0\right)}\left\vert \xx_n^*\left(f_0\right)\right\vert\\
=&t_{j_0}\min_{n\in \supp\left(f_0\right)}\left\vert \xx_n^*\left(f_0+g\right)\right\vert, 
\end{align*}
which implies that $k$ cannot be in a $t_{j_0}$-greedy set of $f_0+g$ of cardinality $d_{j_0}$. Thus, 
\begin{align*}
\left\Vert f\right\Vert\le& \left\Vert f_0\right\Vert\le K\left\Vert f_0+g\right\Vert\le K\left\Vert f+g\right\Vert+K\left\Vert \eta\Ind_{A}\right\Vert
\le K \left\Vert f+g\right\Vert+\mu \left\Vert f\right\Vert,
\end{align*}
which entails that 
\begin{align*}
\left\Vert f\right\Vert\le \frac{K }{1-\mu}\left\Vert f+g\right\Vert\le \left(K+\epsilon\right)\left\Vert f+g\right\Vert.
\end{align*}
By Corollary~\ref{corollarybounded1}, $\langle \XB^*\rangle$ is $\left(K+\epsilon\right)^{-1}$-norming. As $\epsilon$ is arbitrary, it follows that $\langle \XB^*\rangle$ is $K^{-1}$-norming. \\
It remains to consider the case when $d_j=j$ for all $j\in \NN$. In this case we can pick $A=\emptyset$ and $f_0=f$ in the proof above. Thus, again Corollary~\ref{corollarybounded1} gives the desired result.  
\end{proof}
\begin{remark}\rm Note that in  Theorem~\ref{theoremQGequivbidual}, we do not need a Markushevich hypothesis; rather, this property follows at once from our result. The same holds for Theorem~\ref{theoremconvergencenorming}. \end{remark}
In Theorem~\ref{theoremQGequivbidual}, the hypotheses require that for each $f\in \langle \XB\rangle$, we have uniform boundedness of at least a certain sequence of projections $\left(P_{A_j}\left(f\right)\right)_{j\in \NN}$, where $A_j\in \cG\left(f,d_j,t_j\right)$ for each $j\in \NN$. But while convergence of the algorithm and uniform boundedness of the projections are equivalent properties when we are dealing with the TGA \cite[Theorem 4.1]{AABW2021}, \cite[Theorem 1]{Wo2000}, \cite[Theorem 2.1]{Oikhberg2018}, this is not always the case when we consider weaker versions of it. For example, in \cite[Proposition 3.1]{DKSWo2012}, it is shown that given $0<t<1$, for the Haar basis in $L_1\left[0,1\right]$ there is a family $\left\lbrace A_{f,n} \right\rbrace_{\substack{f\in \XX\\n\in \NN}}$ with $A_{f,n}\in \cG\left(f,n,t\right)$ for each $f\in \XX$ and $n\in \NN$ that has the property that $\left(P_{A_{f,n}}\left(f\right)\right)_{n\in\NN}$ converges to $f$ for each $f\in \XX$, but 
\begin{align*}
\left\lbrace \frac{\left\Vert P_{A_{f,n}}\left(f\right)\right\Vert}{\left\Vert f\right\Vert } \right\rbrace_{\substack{f\in \XX\setminus\left\lbrace 0\right\rbrace\\n\in \NN}}
\end{align*}
 is an unbounded family. In our context, we would want a result that holds assuming at most only the described  convergence property of the Haar basis. In our next theorem, we prove this and go further: if we ask only for pointwise boundedness and a sequence $\left(t_j\right)_{j\in \NN}$ instead of fixed $t$, we still get that $\langle \XB^*\rangle$ is norming. To prove our result, we need two auxiliary lemmas: First, a simple lemma that shows how a norming property is inherited from the subsystem $\XB_{>m}$, for any $m\in \NN$. 
\begin{lemma}\label{lemmainherits}Let $\XB$ be a system for $\XX$ and $K\ge 1$. If there is $m\in \NN$ such that $\XB_{>m}^*$ spans a $K^{-1}$-norming subspace of $\XX_{>m}^*$, then there is $K_0\ge K$ such that $\XB^*$ spans a $K_0^{-1}$-norming subspace of $\XX^*$.
\end{lemma}
\begin{proof}
By an inductive argument, we may assume that $m=1$.  Fix $f\in \langle \XB\rangle$  and $\epsilon>0$, and let $f_1:=\xx_1^*\left(f\right)\xx_1$ and $f_2:=f-f_1$. By hypothesis and Lemma~\ref{lemmasimpleeqv4}, there is $m\left(f_2,\epsilon\right)\in \NN_{\ge 2}$ such that for each $g\in\langle \XB_{>m\left(f_2,\epsilon\right)}\rangle$, 
\begin{align*}
\left\Vert f_2\right\Vert\le& \left(K+\epsilon\right) \left\Vert f_2+g\right\Vert.
\end{align*}
Now fix $g \in \langle \XB_{>m\left(f_2,\epsilon\right)}\rangle$. As
\begin{align*}
\left\Vert f+g\right\Vert\ge& \frac{\left\vert \xx_1^*\left(f+g\right)\right\vert}{\left\Vert \xx_1^*\right\Vert}=\frac{\left\Vert f_1\right\Vert}{\left\Vert \xx_1^*\right\Vert\left\Vert \xx_1\right\Vert},
\end{align*}
we have 
\begin{align*}
\left\Vert f\right\Vert\le & \left\Vert f_1\right\Vert+  \left(K+\epsilon\right)\left\Vert f_2+g\right\Vert\le \left(K+1+\epsilon\right)\left\Vert f_1\right\Vert+\left(K+\epsilon\right) \left\Vert f+g\right\Vert\\
\le& \left(K+\epsilon+\left(K+1+\epsilon\right)\left\Vert \xx_1^*\right\Vert\left\Vert \xx_1\right\Vert\right)\left\Vert f+g\right\Vert.
\end{align*}
Thus, Lemma~\ref{lemmasimpleeqv4} gives that $\langle \XB^*\rangle$ is $\left(K+\left(K+1\right)\left\Vert \xx_1^*\right\Vert\left\Vert \xx_1\right\Vert\right)^{-1}$-norming. 
\end{proof}
Our second auxiliary lemma will allow us to pick the vectors in the proof of our theorem with more flexibility. 
\begin{lemma}\label{lemmanormingequivX>m}
Let $\XB$ be a system for $\XX$ with bounded dual system $\XB^*$. If $\XB^*$ does not span a norming subspace of $\XX^*$, for every $\epsilon, r>0$ and $m\in \NN$, there is $f\in \XX$ with the following properties: 
\begin{enumerate}[\rm (i)]
\item \label{rightshifti} $f\in  r S_{\langle \XB_{>m}\rangle}$. 
\item \label{smallinftynorm} $\left\Vert f\right\Vert_{\infty}<\epsilon$. 
\item \label{approximable} For each $n>m$ and each $\delta>0$ there is $g\in \langle \XB_{>n}\rangle$ such that  
\begin{align*}
&\left\Vert g\right\Vert_{\infty}<\delta &&\text{and}&&&\left\Vert f+g\right\Vert<\epsilon. 
\end{align*}
\end{enumerate}
\end{lemma}
\begin{proof}
By Lemma~\ref{lemmainherits}, $\langle \XB_{>m}^*\rangle$ does not span a norming subspace of $\XX_{>m}^*$. Thus, by Corollary~\ref{corollarybounded1}, for each $k\in \NN$ there exists $f_k\in \langle \XB_{>m}\rangle$ with the property that, for each $j>m$, there is $g_{k,j}\in \langle \XB_{>j}\rangle$ with $\left\Vert g_{k,j}\right\Vert_{\infty}<\frac{1}{j}$ such that $\left\Vert f_k\right\Vert>k\left\Vert f_k+g_{k,j}\right\Vert$. For every $k,j\in \NN$, let 
\begin{align*}
&F_k:= r \frac{f_k}{\left\Vert f_k\right\Vert}&&\text{and}&&& G_{k,j}:=r \frac{g_{k,j}}{\left\Vert f_k\right\Vert}. 
\end{align*}
For each $k,j\in \NN$ with $j>\supp\left(F_k\right)$, 
\begin{align*}
&\left\Vert F_k\right\Vert=r>k\left\Vert F_k+G_{k,j}\right\Vert;\\
&\left\Vert F_k\right\Vert_{\infty}\le \left\Vert F_k+G_{k,j}\right\Vert_{\infty}\le \alpha_2 \left\Vert F_k+G_{k,j}\right\Vert< \frac{\alpha_2 r}{k};\\
&\left\Vert G_{k,j}\right\Vert_{\infty}<\frac{r}{\left\Vert f_k\right\Vert j}. 
\end{align*}
Hence, if we choose $k$ so that 
\begin{align*}
\frac{r}{k}\max\left\lbrace 1,\alpha_2\right\rbrace <\epsilon, 
\end{align*}
it follows easily that $f:=F_k$ has the desired properties. 
\end{proof}
We close this section with our second theorem. 
\begin{theorem}\label{theoremconvergencenorming}Let $\XB$ be a system for $\XX$ with dual system $\XB^*$, $\left(d_j\right)_{j\in \NN}$ an unbounded sequence of positive integers, and $\left(t_j\right)_{j\in\NN}$ a scalar sequence with $0<t_j<1$ for all $j\in \NN$. Suppose that for every $f\in \XX$ there is a sequence $\left(A_{f,j}\right)_{j\in \NN}$ with $A_{f,j}\in \cG\left(f,d_{j},t_{j}\right)$ for all $j\in \NN$ such that $\left(P_{A_{f,j}}\left(f\right)\right)_{j\in \NN}$ is bounded. Then $\XB^*$  spans a norming subspace of $\XX^*$. 
\end{theorem}
To simplify our notation in the proof, given $m\in \NN$ and $\epsilon, r>0$, we will say that $f\in \XX$ has Property $Q\left(r,m,\epsilon\right)$ if \ref{rightshifti}, \ref{smallinftynorm} and \ref{approximable} of Lemma~\ref{lemmanormingequivX>m} hold. 
\begin{proof}
The boundedness of $\XB^*$ is proved by the same argument as in Theorem~\ref{theoremQGequivbidual}. \\
Now suppose, for a contradiction, that $\XB^*$ does not span a norming subspace. Let  $f_0=g_0:=0$, $A_0:=\emptyset$ and $m_1=1$, choose $0<\epsilon_1< \frac{1}{2}$ and apply Lemma~\ref{lemmanormingequivX>m} to obtain $f_1$ with Property $Q\left(2, m_1,\epsilon_1\right)$. Choose $j_1$ so that $d_{j_1}>\left\vert \supp\left(f_1\right)\right\vert$, and then $\supp\left(f_1\right)<A_1\in \NN^{<\infty}$ and $0<\eta_1<1$ so that the following hold: 
\begin{align*}
&d_{j_1}=\left\vert \supp\left(f_1\right)\cup A_1\right\vert;&& \eta_1< \min_{n\in \supp\left(f_1\right)}\left\vert \xx_n^*\left(f_1\right)\right\vert; &&&\eta_1\left\Vert\Ind_{A_1}\right\Vert< \epsilon_1.
\end{align*}
Now pick $g_1\in \langle \XB_{>\max\left(A_1\right)}\rangle$ so that
\begin{align*}
&\left\Vert g_1\right\Vert_{\infty}<\eta_1 t_{j_1} &&\text{and}&&& \left\Vert f_1+g_1\right\Vert< \epsilon_1, 
\end{align*}
let $m_2:=\max\left(\supp\left(g_1\right)\right)$, and choose 
\begin{align*}
0<\epsilon_2< \frac{1}{2}\min\left\lbrace \epsilon_1, \min_{n\in \supp\left(g_1\right)}\left\vert \xx_n^*\left(g_1\right)\right\vert \right\rbrace. 
\end{align*}
Apply Lemma~\ref{lemmanormingequivX>m} to get $f_2$ with Property $Q\left(4, m_2, \epsilon_2\right)$, then choose $j_2>j_1$ so that
\begin{align*}
d_{j_2}>\left\vert \supp\left(f_2\right) \cup \supp\left(f_1\right)\cup A_1 \cup \supp\left(g_1\right) \right\vert,
\end{align*}
choose $A_2>\supp\left(f_2\right)$ so that 
\begin{align*}
d_{j_2}=&\left\vert \supp\left(f_2\right)\cup A_2 \cup \supp\left(f_1\right)\cup A_1 \cup \supp\left(g_1\right) \right\vert,
\end{align*}
and continue as before so that, in a recursive manner, we construct sequences $\left(f_k\right)_{k\in \NN}, \left(g_k\right)_{k\in \NN}$, $\left(\epsilon_k\right)_{k\in \NN}, \left(\eta_k\right)_{k\in \NN}$, and $\left(m_k\right)_{k\in \NN}, \left(j_k\right)_{k\in \NN}$ such that for each $k\in \NN$, the following hold:
\begin{align*}
&\supp\left(f_k\right)<A_k
<\supp\left(g_k\right)<\supp\left(f_{k+1}\right);\\
&\left\vert \supp\left(f_k\right)\cup A_k \bigcup_{1\le l\le k-1}\supp\left(f_l\right)\cup A_l\cup \supp\left(g_l\right) \right\vert=d_{j_k};\\
&j_k<j_{k+1}\qquad\text{and} \qquad m_k<m_{k+1};\\
&\eta_k \left\Vert \Ind_{A_k}\right\Vert< \epsilon_k\qquad\text{and}\qquad \eta_k< \min_{n\in \supp\left(f_k\right)}\left\vert \xx_n^*\left(f_k\right)\right\vert;\\
&\left\Vert g_k\right\Vert_{\infty}<\eta_k t_{j_k} \qquad \text{and}\qquad\left\Vert f_k+g_k\right\Vert< \epsilon_k;\\
&\epsilon_{k+1}< \frac{1}{2}\min\left\lbrace \epsilon_k, \min_{n\in \supp\left(g_k\right)}\left\vert \xx_n^*\left(g_k\right)\right\vert\right\rbrace;\\
&f_k \qquad\text{has Property}\qquad Q\left(2^k, m_k,\epsilon_k\right).
\end{align*}
For each $k\in \NN$, define  
\begin{align*}
&B_k:=  \supp\left(f_k\right)\cup A_k \bigcup_{1\le l\le k-1}\supp\left(f_l\right)\cup A_l\cup \supp\left(g_l\right);\\
&D_k:=B_k\cup \supp\left(g_k\right),
\end{align*}
and let 
\begin{align*}
&f_0:=\sum_{k=1}^{\infty}\left(f_k+\eta_k\Ind_{A_k}+g_k\right). 
\end{align*}
Note that $f_0$ is well-defined, with 
\begin{align*}
&P_{D_k}\left(f_0\right)\xrightarrow[k\to \infty]{}f_0. 
\end{align*}
For each $k\in\NN$,
\begin{align*}
\left\Vert P_{D_k}\left(f_0\right)-P_{B_k}\left(f_0\right)\right\Vert=\left\Vert g_k\right\Vert\ge \left\Vert f_k\right\Vert -\left\Vert f_k+g_k\right\Vert > 2^{k}-\epsilon_k>2^{k-1}.
\end{align*}
Hence, $\left(P_{B_k}\left(f_0\right)\right)_{k\in \NN}$ is unbounded. Thus, to derive a contradiction and finish the proof, it suffices to show that for each $k\in \NN$, the only $t_{j_k}$-greedy set of $f_0$ of cardinality $d_{j_k}$ is $B_k$. To this end, first note that for each $k\in \NN$,
\begin{align*}
&\left\Vert f_{k+1}\right\Vert_{\infty}<\min_{n\in \supp\left(g_k\right)}\left\vert \xx_n^*\left(g_k\right)\right\vert\le \left\Vert g_{k}\right\Vert_{\infty}<t_{j_k}\eta_k\le  t_{j_k} \min_{n\in \supp\left(f_k\right)}\left\vert \xx_n^*\left(f_k\right)\right\vert.
\end{align*}
Now fix $k\in \NN$ and $n\in \supp\left(f_0\right)\setminus B_k$. Note that $n>B_k$, so either $n\in \supp\left(g_k\right)$ or there is $l>k$ such that $n\in \supp\left(f_l\right)\cup A_l\cup \supp\left(g_l\right)$. In either case, it follows from our construction that  
\begin{align*}
&\left\vert \xx_n^*\left(f_0\right)\right\vert< t_{j_k}\min_{i\in B_k}\left\vert \xx_i^*\left(f_0\right)\right\vert. 
\end{align*}
As $\left\vert B_k\right\vert =d_{j_k}$, this entails that $n$ cannot be in a $t_{j_k}$-greedy set of $f_0$ of cardinality $d_{j_k}$. Since $\cG\left(f_0,d_{j_k},t_{j,k}\right)\not=\emptyset$, the only remaining possibility is that $\cG\left(f_0,d_{j_k},t_{j,k}\right)=\left\lbrace B_k\right\rbrace$, and the proof is complete. 
\end{proof}

 \section{Some remarks on other greedy-like properties. }\label{section closing remarks}
A key property in the proof of Theorem~\ref{theoremconvergencenorming} is the boundedness of some (weak) greedy sums, but one may wonder whether some greedy-like properties from the literature that do not involve such boundedness may still share the norming property of our theorems. It turns out that the answer is negative: First, a semi-greedy system may have a dual system that is not total \cite[Example 4.5]{BL2020}, so in particular it does not span a norming subspace. On the other hand, semi-greedy Markushevich bases are almost-greedy \cite[Corollary 4.3]{BL2020}, so there is uniform boundedness and convergence of the greedy sums. The strongest greedy-like property of  interest in the field left to be considered is bidemocracy:  as in the case of semi-greedy systems, bidemocratic systems need not be Markushevich bases \cite[Theorem 3.6, Corollary 3.7]{AABBL2023}. However, unlike the case of semi-greedy systems, a Markushevich hypothesis does not guarantee a norming dual basis. To see this, let us recall some definitions from the literature (see for example \cite{AABW2021}, \cite{DKKT2003}). 
\begin{definition}\label{definitionfundamental}Given a bounded system $\XB$ for $\XX$, the \emph{upper superdemocracy function} $\varphi_u^{\varepsilon}: \NN\rightarrow \RR$ of $\XB$ is defined as
\begin{align*}
\varphi_u^{\varepsilon}\left(n\right)=\varphi_u^{\varepsilon}\left(n\right)\left[\XB,\XX\right]:&\sup_{\substack{A\in \NN^{\le n}\\\varepsilon\in S_{\FF}^{n}}}\left\Vert\Ind_{\varepsilon, A}\right\Vert. 
\end{align*}
\end{definition}
\begin{definition}Let $f:\NN\rightarrow \RR_{>0}$ be a function. We say that $f$ has the \emph{upper regularity property} (URP) if there is $b\in \NN$ such that
\begin{align*}
&f\left(bm\right)\le \frac{b}{2} f\left(m\right)&&\forall m\in \NN, 
\end{align*}
and has the \emph{lower regularity property} (LRP) if there is $b\in \NN$ such that 
\begin{align*}
&2f\left(m\right)\le f\left(bm\right)&&\forall m\in \NN. 
\end{align*}
\end{definition}
We also need a result based on the DKK method for constructing bases, developed first by Dilworth, Kalton and Kutzarova in \cite{DKK2003}. The following lemma is obtained with one of the extensions of the method.  

\begin{lemma}\label{lemmaDKK}Let $\XB$ be a bounded Markushevich basis for $\XX$ such that $\XB^*$ bounded. Suppose that $\XX$ has a complemented subspace $\ZZ$ with a subsymmetric basis $\ZB$ whose upper superdemocracy function has both the LRP and the URP. Then $\XX$ has a bidemocratic Markushevich basis $\YB$ with a block basis equivalent to $\XB$. 
\end{lemma}
\begin{proof}
An application of \cite[Lemma 3.2, Proposition 3.3]{AAB2023s} to $\XB$ and $\ZZ$ gives a a space $\YY$ that is isomorphic to $\ZZ\oplus \XX\approx \XX$ with a bounded Markushevich basis $\YB$ for $\YY$ that has a block basis equivalent to $\XB$. A combination of \cite[Proposition 3.3 (xi)]{AAB2023s} and \cite[Proposition 9.4]{AABW2021} gives that $\YB$ is truncation quasi-greedy and democratic, whereas the same combination together with \cite[Proposition 3.3 (ii)]{AAB2023s} entail that $\YB^*$ also has those two properties. Therefore, by \cite[Corollary 2.6]{AABBL2023}, $\YB$ is bidemocratic. 
\end{proof}

In order to apply Lemma~\ref{lemmaDKK} and prove the existence of bidemocratic Markushevich bases with nonnorming dual bases, we need some suitable bases and spaces from the literature. We can find such examples in \cite{AABW2021} and \cite{DL1972} (see also \cite{HR2024}).

\begin{proposition}\label{propositionbidem} (cf. \cite[Theorem 2.5, Remark 2.6]{HR2024}).  There is a Banach space $\XX$ with a bidemocratic Markushevich basis $\XB$ whose dual basis $\XB^*$ does not span a norming subspace. Moreover, if $\YY$ is separable Banach space with the following properties:  
\begin{enumerate}[\rm (A)]
\item \label{quasi-ref} $\YY$ is not quasi-reflexive.  
\item \label{separabledual} $\YY^*$ is separable. 
\item \label{comp} $\YY$ has a complemented subspace $\ZZ$ that has a subsymmetric basis $\ZB$ whose upper democracy function has the LRP and the URP.
\end{enumerate}
Then $\YY$ has a bidemocratic Markushevich basis $\YB$ such that $\YB^*$ does not span a norming subspace. 
\end{proposition}
\begin{proof}
\cite[Proposition 2.11]{AABW2021} gives a Banach space $\XX$ that has a complemented copy of $\ell_2$ and a Markushevich basis $\XB$ whose dual basis $\XB^*$ is bounded and does not span a norming subspace. By Lemma~\ref{lemmaDKK}, $\XX$ has a  bidemocratic Markushevich basis $\YB$ with a block basis equivalent to $\XB$. By Remark~\ref{remarkblockbasis}, $\YB^*$ does not span a norming subspace. \\
Now suppose $\YY$ is separable and \ref{quasi-ref}, \ref{separabledual} and \ref{comp}  hold. By \cite[Theorem]{DL1972}, there is a total sequence  $\left(y_n^*\right)_{n\in \NN}\subset S_{\XX^*}$ whose linear span is not norming.  Let $\left(y_n\right)_{n\in \NN}\subset \YY$ be a sequence whose linear span is dense in $\YY$. Applying the procedure of \cite[Theorem 1.f.4]{LT1977} to $\left(y_n\right)_{n\in \NN}$ and $\left(y_n^*\right)_{n\in \NN}$ one obtains a bounded Markushevich basis $\YB$ for $\YY$ such that $\YB^*$ is bounded and $\langle \YB^*\rangle$ is not norming, and the proof is completed by the same argument as before. 
\end{proof}
While a bidemocratic minimal system - even if it is a Markushevich basis - may fail to be equivalent to its bidual system, bidemocracy is still passed to the latter \cite[Lemma 5.6]{AABW2021}, so one may ask whether some greedy-like properties weaker than bidemocracy also pass to the bidual system.  It was proved in \cite[Corollary 2.6]{AABBL2023} that a minimal system $\XB$ is bidemocratic if and only if both $\XB$ and $\XB^*$ are truncation quasi-greedy and (partially) democratic, or conservative. We will show that truncation quasi-greediness passes to the bidual system and, if $\XB$ is truncation quasi-greedy, then democracy, superdemocracy, partial democracy and conservativeness all pass to the bidual system as well. To this end, we will use the characterization of truncation quasi-greedy systems from \cite[Theorem 1.5]{AAB2024}, for which we need the following definitions. 
 \begin{definition}\label{definitionoscillation}\cite{DOSZ2009}, \cite{AAB2024} Let $\XB$ be a basis of a Banach space $\XX$, $f\in \XX$ and $A\in \NN^{<\infty}$. The \emph{oscillation} of $f$ on $A$ is given by
\begin{align*}
o\left(f,A\right):=&\max_{\substack{n\in A,k\in A\cap\supp\left(f\right)}}\frac{\left\vert\xx_n^*\left(f\right)\right\vert}{\left\vert\xx_k^*\left(f\right)\right\vert}
\end{align*}
with the convention that $o\left(f,A\right)=0$ if $A\cap\supp\left(f\right)=\emptyset$ . 
\end{definition}
The above definition was used in \cite{DOSZ2009} to define $(D,d)$-bounded-oscillation unconditional systems for some $D, d\ge 1$. For our purposes, it is more convenient to define bounded-oscillation unconditional systems using the equivalence proved in \cite{AAB2024}. 
\begin{definition}\label{definitionBOU}\cite{AAB2024} (cf. \cite{DOSZ2009}). Let $\XB$ be a system $\XX$. We say that $\XB$ is bounded-oscillation unconditional if there is a function $\Phi=\Phi\left[\XB,\XX\right]: \left(0,1\right]\rightarrow \RR$ such that 
\begin{align*}
&\left\Vert P_A\left(f\right)\right\Vert\le \Phi\left(t\right)\left\Vert f\right\Vert&&\forall f\in \XX, A\in \NN^{<\infty} : \osc\left(f,A\right)\le\frac{1}{t}.
\end{align*}
\end{definition}
It was proved in \cite[Theorem 1.5]{AAB2024} that a system is truncation quasi-greedy if and only if it is bounded-oscillation unconditional. We will use this equivalence in the proof of the following lemma -with which we close this note. 
\begin{lemma}\label{lemmaTQG}Let $\XB$ be a truncation quasi-greedy system for $\XX$. For every $f\in \langle \XB\rangle$, $0<t\le 1$ and each $A\subset \supp\left(f\right)$ with $\osc\left(f,A\right)\le \frac{1}{t}$, we have
\begin{align}
\left\Vert P_A\left(f\right)\right\Vert_{\XX}\le&\Phi\left[\XB,\XX\right]\left(t\right)\left\Vert j_{\XX}\left(f\right)\big|_{\overline{\langle \XB^*\rangle}}\right\Vert_{\overline{\langle \XB^*\rangle}^*}\label{bou0}
\end{align}
Thus, $\XB^{**}$ is a truncation quasi-greedy system for the closure of its span in $\overline{\langle \XB^*\rangle}^{*}$, with 
\begin{align}
\Phi\left[\XB^{**},\overline{\langle \XB^{**}\rangle}\right]\left(t\right)\le \Phi\left[\XB,\XX\right]\left(t\right)&&\forall 0<t\le 1.\label{bou1}
\end{align}
In particular, if $\XB$ is also partially democratic, conservative, democratic or superdemocratic, so is $\XB^{**}$.
\end{lemma}
\begin{proof}
We proceed as in the proof of Lemma~\ref{lemmasimpleeqv4}: Given $f$ and $A$ as in the statement, define $g^*: \langle P_A\left(f\right)\rangle\oplus \langle \xx_n: n\not\in A\rangle$ by 
\begin{align*}
&g^*\left(a P_A\left(f\right)+ g\right)=a\left\Vert P_A\left(f\right)\right\Vert_{\XX}&&\forall g\in \langle \xx_n: n\not\in A\rangle. 
\end{align*}
Since $\osc\left(f,A\right)\le t^{-1}$, it follows that $\left\Vert g^*\right\Vert\le \Phi\left[\XB,\XX\right]\left(t\right)$. Let $f^*$ be a norm-preserving extension of $g^*$ to $\XX$. As in the proof of Lemma~\ref{lemmasimpleeqv4}, we get that $f^*=\sum_{n\in A}f^*\left(\xx_n\right)\xx_n^*$. Since
\begin{align*}
\Phi\left[\XB,\XX\right]\left(t\right)\left\Vert j_{\XX}\left(f\right)\big|_{\overline{\langle \XB^*\rangle}}\right\Vert_{\overline{\langle \XB^*\rangle}^*}\ge \left\vert j_{\XX}\left(f\right)\big|_{\overline{\langle \XB^*\rangle}}\left(f^*\right)\right\vert=\left\Vert P_A\left(f\right)\right\Vert_{\XX}, 
\end{align*}
we have obtained  \eqref{bou0}, and then \eqref{bou1} follows by density. The rest of the statements follow from the fact that for every $B\in \NN^{<\infty}$ and every $\varepsilon\in S_{\FF}^{B}$, by \eqref{bou0} we have
\begin{align*}
\left\Vert\Ind_{\varepsilon, B}\left[\XB^{**},\overline{\langle \XB^{**}\rangle}\right]\right\Vert_{\overline{\langle \XB^*\rangle}^{*}}\le& \left\Vert \Ind_{\varepsilon, B}\right\Vert\left[\XB,\XX\right]\\
\le& \Phi\left[\XB,\XX\right]\left(1\right) \left\Vert\Ind_{\varepsilon, B}\left[\XB^{**},\overline{\langle \XB^{**}\rangle}\right]\right\Vert_{\overline{\langle \XB^*\rangle}^{*}}.
\end{align*}
\end{proof}

\section*{Declarations}

\subsection*{Funding}

The author was supported by  the Grants ANPCyT PICT 2018-04104 and CONICET PIP 11220200101609CO.

\subsection*{Conflict of Interest/Competing Interest}

The author has no conflict of interest or competing interest.

\end{document}